\documentclass
[12pt]{amsart}

\usepackage{a4}
\usepackage{amsthm,amsfonts,amssymb,latexsym}
\usepackage{color}
\usepackage{tikz}
\usepackage{tikz-qtree}

\def\ab{\{a,b\}}
\def\AB{\{A,B\}}
\def\lt{\triangleleft}
\def\rt{\triangleright}
\def\rlt{\bowtie}
\def\prf{\begin{proof}}
\def\prfend{\end{proof}}

\def\LS{\mathcal{LS}}
\def\A{\mathbb{A}}
\def\S{\mathcal{S}}
\def\R{\mathbb{R}}
\def\sgn{{\rm sgn}}
\def\fsgn{{\rm SGN}}
\def\tr{{\rm tr}}

\newtheorem{lem}{Lemma}[section]
\newtheorem{prp}{Proposition}[section]
\newtheorem{thm}{Theorem}[section]
\newtheorem{cor}{Corollary}[section]
\newtheorem{definition}{Definition}[section]

\title[Complexity of Substitutive Sequences]{
Calculation of the Complexities of Substitutive Sequences Over a Binary Alphabet
}

\author{Bo TAN}
\address{School of Mathematics and Statistics, Huazhong University of Science and Technology, 430074 Wuhan, P.~R.~China}
\email{tanbo@hust.edu.cn}

\author{Zhi-Xiong WEN}
\address{School of Mathematics and Statistics, Huazhong University of Science and Technology, 430074 Wuhan, P.~R.~China}
\email{zhi-xiong.wen@hust.edu.cn}

\author{Yiping ZHANG}
\address{ School of Mathematics and statistics, Wuhan University, 430072 Wuhan, P.~R.~China}
\email{ypzhang@whu.edu.cn}

\thanks{Research supported by NSFC. 11171123 \& 11222111.}

\date{}

\begin{document}

\maketitle

\begin{abstract}

We consider the complexities of   substitutive sequences  over a  binary alphabet. By studying various types of special words, we show that,  knowing some initial values, its complexity   can be completely formulated via a recurrence formula determined by the characteristic polynomial.



\end{abstract}

\medskip



\section{Introduction }

The study of substitutions over a finite alphabet plays  important
roles in many fields  such as finite automata, symbolic dynamics, formal languages,
number theory, fractal geometry {\it etc.} It has various applications
to quasi-crystals, computational complexity, information
theory\ldots (see \cite{all,AI,Co,Lo,Lo2} and the references
therein).  In addition, substitutions are also fundamental objects in
combinatorial group theory \cite{
Ly,MKS}.

\smallskip

Given an infinite sequence $\xi=\xi_1\xi_2\xi_3\cdots$ ($\xi_i\in \A $) over some finite alphabet
$\A $, we denote by ${\mathcal
  L}_n(\xi)$ the set $\{\xi_{i}\cdots \xi_{i+n-1}\bigm| i\ge1\}$ of
factors of~$\xi$ of length~$n\ (n\ge 1)$, and by convention ${\mathcal
  L}_0(\xi)$ is the singleton consisting of the empty word $\varepsilon$. The set ${\mathcal
  L}(\xi)=\cup_{n\ge 0}{\mathcal L}_n(\xi)$ is then called the language of
$\xi$, and the function $p_\xi(n):=\#{\mathcal L}_n(\xi)$ the
complexity of $\xi$, here and hereafter $\#$ denotes the cardinality of a finite set.

Let $\A ^*$ be the free monoid generated by $\A $ (with   $\varepsilon$ as the neutral  element). A morphism $\sigma:
\A ^* \to \A ^*$ is called a substitution.
We  deal with only the non-erasing substitutions (the image of any letter in $\A $ is not the empty word), whence the substitution can be extended naturally to $\A ^\mathbb{N}$, the set of infinite sequences over $\A $.
 Denote by
$\xi_{\sigma}$ any one of the fixed points of $\sigma$ (that
is~$\sigma(\xi_{\sigma})=\xi_{\sigma}$), if it exists.

\smallskip

The study of the complexity of $\xi_{\sigma}$ (also called the complexity of ${\sigma}$) has a long history. In general, it is very difficult to find out the explicit formula for $p_{\xi}(n)$ for a given ${\sigma}$; only some calculations for specific classes of substitutions can be found in the literature. Here are some known results :

\begin{itemize}
\item $p_\xi(n)\le n$ for some $n$ if and only if $\xi$ is ultimately
  periodic, and in this case the complexity is bounded~\cite{MH};
\item A sequence $\xi$  of complexity
  $p_{\xi}(n)=n+1$ is called Sturmian. There are many equivalent
  characterizations and interesting properties of Sturmian sequences (see, e.g.
  \cite{Lo,Se,WWli});
\item Rote \cite{Ro} constructed a class of sequences with complexity
  $2n$ by using graphs;
\item Moss\'e \cite{Mo} studied the case of $q$-automata
  (which correspond to substitutions of constant length). A method to compute $p(n)$ with
  linear recurrence formula was given under some technical
  conditions;
\item Over a ternary alphabet, a class of Tribonacci type
  substitutions with complexity $2n+1$ was introduced by Arnoux and
  Rauzy \cite{AR}.
 An example of substitution (Triplex Substitution)  with complexity $3n$ is presented by the authors  \cite{twz1}.

 \item For a fixed point of some substitution, the complexity can only be of the following five different asymptotic forms: $\Theta(1), \Theta(n), \Theta(n\log\log n), \Theta(n\log n)$ or $\Theta(n^2)$, where $\Theta(g(n))$ means a function $f(n)$ satisfying $0<\liminf \frac{f(n)}{g(n)}\le
     \limsup \frac{f(n)}{g(n)}<\infty$ \cite{Pa}.

     \item For a survey and more general computation of factor complexity of word (on a alphabet of cardinality more than 2), we suggest to see \cite{CaN,Fr}.
\end{itemize}


In this paper, we consider general substitutions $\sigma$ over a binary alphabet. Using Moss\'e's theory of identifiability (\cite{Mo}) and by studying various types of special words (\cite{Ca1,CaN}), we show that the complexity $p(n)$ can be completely formulated knowing some initial values, and a recurrence formula is given.

\section{Notations and Preliminary }

We fix the binary alphabet ${{\A }}=\{a, b\}$ consisting of two letters  $a$ and $b$. Let
${\A }^*$ be the free monoid generated by ${\mathbb
  A}$ (with the empty word $\varepsilon$ as the neutral element),
 and ${\A }^\mathbb{N}$ be the set of all infinite sequences (also called infinite words) over $\A$. 

If $w\in {\A }^{*}$, we denote by $|w|$ its length and by
$|w|_{a}$ (resp. $|w|_{b}$) the number of occurrences of
the letter~$a$ (resp.~$b$) in $w$. The abelian Parikh vector of $w$ is then defined to be the column vector
 $L(w)=(|w|_a, |w|_b)^t\in \mathbb{N}^2$.

\smallskip

A word $v$ is a factor of a word $w$ (written as $v\in w$) if
there exist $u,\, u'\in {\A }^*$, such that $w = uvu'$ . It is sometimes convenient to use the notation
``$\circledast$" to stand for some word 
which we don't care so much. Thus $v$ is a factor of a word $w$ if and only if $w=\circledast v  \circledast $ (remark that even within a formula, $\circledast$'s may represent different words).  We say
that $v$ is a prefix  (resp. a suffix) of $w$ if $w=v  \circledast $ (resp.~$w=\circledast  v$), and then we
write $v\triangleleft w$  (resp. $v\triangleright w$).
Two words $v$ and $w$ are said to be comparable, written $v\bowtie w$, if  either $v\rt w$ or $w\rt v$.
 The notions of
factor and prefix extend to infinite words in a natural way.

It is also convenient to put, e.g. $\A^*v:=\{xv; x\in \A^*\}$,
$\A^*v\A^*:=\{xvy; x,y\in \A^*\}$, {\it etc}. Thus $w\in \A^*v\A^*\Leftrightarrow v\in w$; $w\in v\A^*\Leftrightarrow v\lhd w$, and so on.

When $\xi=w_1\cdots w_m\cdots\in \A^*\cup \A^\mathbb{N}~ (w_i\in\A)$, we   also write $\xi|_1=w_1, \cdots, \xi|_m=w_m, \cdots$, and $\xi[i,j]=w_iw_{i+1}\cdots w_j (i\le j)$.

\smallskip

As already defined, a substitution $\sigma$ over~$\A$ is a morphism
$\sigma$ of $\A^*$. 
The matrix $M=
(L(\sigma(a)), L(\sigma(b)))$ is called the incidence
matrix of~$\sigma$. The characteristic polynomial $\lambda^2-\tr(M)\lambda+\det(M)$ of $M$ is also called the characteristic polynomial of $\sigma$.

If $\sigma(a)$ and $\sigma(b)$ have distinct first letters, we say that the substitution $\sigma$ is marked, and if moreover $\sigma(a)=a\circledast $ and $\sigma(b)=b\circledast $, we say that $\sigma$ is well-marked. It is easy to see that  $\sigma^2$ is well-marked if $\sigma$ is marked.

In this paper, all substitutions are assumed to be non-erasing, that is, the image of each letter is not empty. Whence, the substitution can be extended naturally to $\A^\mathbb N$.
An infinite word  $\xi=\xi_1\xi_2\cdots$ is a fixed point
of $\sigma$ if $\sigma(\xi)=\xi$.

\smallskip

Hereafter, we suppose  that the substitution $\sigma$ is primitive (i.e. its incidence matrix $M$
is primitive: $M^n$ possesses positive coordinates for some positive integer $n$).
 The following easy facts for a primitive substitution $\sigma$  are well known:

 \begin{enumerate}

   \item the fixed point of $\sigma$ is recurrent, that is, every factor will occur for infinitely many times; and all the fixed points of  $\sigma$ have the same language;
   \item  a substitution $\sigma$ and its powers $\sigma^n$ ($n\ge1$) have the same fixed points, and thus have the same language;
   \item if one substitution is a composition of an inner automorphism (of the free group) with another substitution, then the two substitutions have the same language.
 \end{enumerate}
We suppose also that the fixed point $\xi$ of $\sigma$ is not (ultimately) periodic;
the periodic case are characterized completely by S\'e\'ebold \cite{Se2}. In
particular, whence $\{\sigma(a), \sigma(b)\}$ is a code,
and thus $\sigma$ is marked up to an inner automorphism (see \cite{Lo}).
For the sake of calculation of the complexity of a non-periodic primitive substitution, we may  further suppose, without   loss of generality, that  the substitution is well-marked.

\smallskip

The notion of ``special words" is a powerful tool  for calculating the
complexity.  See \cite{Ca1,CaN}  and \cite{Ca,Lo,Lo2} for more information.

Let $W$ be a factor of $\xi$. If  $\delta\in {\A }$ such that $W\delta$ is a factor of $\xi$, then we say that $W\delta$ is a
right extension  of $W$. A word is called a right
special word (special word for short) of $\xi$ if it has more than one
extensions, that is, $Wa\in \xi$ and $Wb\in \xi$. Similarly we define ``left extension" and ``left special
word". It is easy to see that a suffix (resp. prefix) of a special (resp. left special) word is also  special (resp. left special).

Let $\mathcal{S}_n$ (resp. $\mathcal{LS}_n$) be the set of special words (resp. left special words) of length $n$ of $\xi$. Put $\mathcal{S}=\cup_{n\ge0}\mathcal{S}_n$ (resp. $\mathcal{LS}=\cup\mathcal{LS}_n$). It is easy to see that
 $$s(n):=\#\mathcal{S}_n=\#\mathcal{LS}_n=\Delta p(n+1)(:=p(n+1)-p(n)).$$ Hence the study of $p(n)$ is almost equivalent to the study of $s(n)$.

\medskip

\subsection{The word  $W_0$ and the letters $\delta_a,\delta_b$} \mbox{}

\smallskip

 Write $A=\sigma(a), B=\sigma(b)$, and denote $\AB^*$
  the set of words obtained by a finite concatenation of the words
   $A$ and $B$. Put, as before, e.g.
 $\AB^* A:=\{VA; V\in\AB^*\}.$ 
Remark that since $\sigma$ is non-periodic,  $\{A, B\}$ is a code  and
$\AB^*$ is a disjoint union of $\AB^*A$ and $\AB^*B$.

 Since $\sigma$ is non-periodic, the left-infinite words $A^{\infty}(=\cdots AA\cdots A)$ and $B^{\infty}$ are different. Let $W_0$ be the longest common suffix of $A^{\infty}$ and $B^{\infty}$ (see also  \cite{tan}). Remark that $W_0$ is possibly empty.

  The following lemma is a direct  consequence of Fine-Wilf theorem \cite{Fogg}.

\begin{lem}\label{l2.1}
$|W_0|\le |A|+|B|-2$.
\end{lem}


By the definition of $W_0$, for some 
$\delta_a, \delta_b \in \ab$
with $\{\delta_a,\delta_b \}=\ab$,
 \begin{equation}\label{alpha}
 A^{\infty}=\circledast \delta_a W_0 {\textrm \ \ and\  \  } B^{\infty}=\circledast \delta_b W_0.
 \end{equation}


 Formula (\ref{alpha}) shows that there exist  $m\ge 0$ and $A'\rt A$ ($|A'|<|A|$)  such that
\begin{equation}\label{e2.2}
W_0=A'A^m ,  {\textrm \ \ and\  \  }  \delta_a A' \rt A,
\end{equation}
and similarly
$$W_0=B'B^k ,  {\textrm \ \ and\  \  }  \delta_b B' \rt B.$$

The following lemma is essentially due to \cite{tan}.
 \begin{lem}\label{w0}
  (1) For   $W\in\AB^*$,  we have $W_0\rlt W$. 
  Furthermore,

 (2)  If $W\in \AB^*A$ (resp. $ \AB^*B$ ) and $|W|> |W_0|$, then
  $ \delta_a W_0\rt W $ (resp. $  \delta_b W_0\rt W $), where $\delta_a$ and $\delta_b$ are defined in (\ref{alpha}).

(3) Let $W\in\AB^*$.  If $\delta_aW_0\rt W$ (resp. $\delta_bW_0\rt W$), then $W\in \AB^*A$ (resp. $\AB^*B$).

In brief, any word in $\AB^*$ is comparable with $W_0$. Amongst them, the word in  $\AB^*A$ is comparable with $\delta_a W_0$ and $\AB^*B$ is comparable with $\delta_b W_0$.
%
%
\end{lem}

\prf 
If $W=A$ or $W=B$, the lemma is obvious.  Suppose   $W\in \AB^*$ such that
$W_0\rlt W$, we claim that $\delta_aW_0\rlt WA$ and $\delta_bW_0\rlt WB$. The two statements can be proven in the same way, and we only show the first one by considering the following two cases:
\smallskip

\noindent Case 1: $W_0\rt W$.  Then  $W_0A\rt WA$, and on the other hand, $\delta_aW_0\rt W_0A$ because both of them are suffixes of $A^\infty$. Hence
$\delta_aW_0\rt WA$.
\smallskip

\noindent Case 2: $W\rt W_0$. Then  $WA\rt W_0A$, while $W_0A$ is a suffix of $A^\infty$, and thus $WA$ is a suffix of $A^\infty$. This yields that $WA\rlt \delta_aW_0$ because both of them are suffixes of $A^\infty$.
\prfend

\begin{cor}\label{cor21}
 Let $W\in\AB^*$. Then $W_0\rt W_0W$, $\delta_aW_0\rt W_0WA$,  $\delta_bW_0\rt W_0WB$. In particular,  $\delta_aW_0\rt W_0A$,  $\delta_bW_0\rt W_0B$.
\end{cor}

\medskip
\subsection{Natural decomposition and identifiability}  \mbox{}

\smallskip

Let $\xi$ be a fixed sequence of $\sigma$. Write $\xi=\xi_1\xi_2\cdots$.
Since $\sigma(\xi)=\xi$, we have the following so called ``natural decomposition" of $\xi$
\begin{equation}\label{e0}
\xi=[\xi_{1}\xi_2\cdots\xi_{n_2-1}][\xi_{n_2}\cdots\xi_{n_3-1}]\cdots
[\xi_{n_{k}}\cdots\xi_{n_{k+1}-1}][\xi_{n_{k+1}}\cdots,
\end{equation}
 where $\xi_k\in \A=\ab $, $\sigma(\xi_k)=\xi_{n_{k}}\cdots\xi_{n_{k+1}-1}\in\AB \ (k\ge1)$, and $n_1(:=1),\cdots, n_k(:=|\sigma(\xi[1,k-1])|+1), \cdots$ are called the ``cutting positions" of $\xi$. We denote
\begin{equation}\label{position1}
E_1=\{n_k; k\ge1\}.
\end{equation}

Now consider the factors of $\xi$. Let $W=\xi_i\xi_{i+1}\cdots\xi_j\in\xi$, then (comparing to (\ref{e0})) for some integers
$ k,l $ $(n_{k-1}< i\le n_k\le n_l\le n_{l+1}-1\le j< n_{l+2})$, we have
$$W=
\xi_i\cdots\xi_{n_k-1}][\xi_{n_k}\cdots\xi_{n_{k+1}-1}]\cdots[\xi_{n_{l}}
\cdots\xi_{n_{l+1}-1}][\xi_{n_{l+1}}\cdots\xi_{j},
$$
\noindent that is, observing the cutting positions of $W$ in $\xi$ we can write out the following natural decomposition of $W$
\begin{equation}\label{e2.9}
W=U\sigma(\xi_k)\cdots\sigma(\xi_l)V=U\sigma(W')V,
\end{equation}
 where
 \begin{align*}
  U&=\xi_i\cdots\xi_{n_k-1}\rt\sigma(\xi_{k-1}),~  |U|<|\sigma(\xi_{k-1})|,\\
  \sigma(\xi_k)&=\xi_{n_k}\cdots\xi_{n_{k+1}-1},\\
  &\vdots \\
  \sigma(\xi_l)&=\xi_{n_l}\cdots\xi_{n_{l+1}-1},\\
  V&=\xi_{n_{l+1}}\cdots\xi_{j}\lt\sigma(\xi_{l+1}), ~ |V|<|\sigma(\xi_{l+1})|,\\
  W'&=\xi_k\cdots\xi_l\in \xi.
 \end{align*}
%
%
%
%
%
%
We say that $W'$ (resp. $\xi_m$ ,\ $k\le m\le l$ ) is the ancestor of $\sigma(W')$ (resp. $\sigma(\xi_m)$). Sometimes, we also call $\xi_{k-1}\xi_k\cdots\xi_l\xi_{l+1}$ the ancestor of $W$.

\medskip
We extend a little more the significance of ``natural decomposition'': if  $W=U\sigma(W_1W''W_2)V$   as in (\ref{e2.9}), we shall also say that $W=U'\sigma(W'')V'$ is a ``natural decomposition'' (where $U'=U\sigma(W_1), V'=\sigma(W_2)V$), and we  write $W= U'[\sigma(W'')]V'.$
Equivalently, the notation $U'[\sigma(W'')]V'$ means that there exist $U'',V''\in\A^*$ such that
 \begin{equation}\label{e210}
 U''W''V''\in\xi, ~U'\rt \sigma(U''), \ \textrm{\ and\ }  V'\lt \sigma(V'') .
 \end{equation}
Intuitively, $U'[\sigma(W'')]V'$ appears in $\xi$ with $``[\;"$ and $``\;]"$ showing the interested natural cutting positions.

We call the decomposition as in (\ref{e2.9}) a strict natural decomposition of $W$. Remark that any natural decomposition can be extended to a strict one, and, in general, the natural decompositions of a factor  are not unique;
and that the fact $U\sigma(W)V\in\xi$ does not always mean   $U[\sigma(W)]V$  !

\medskip

From the theory of identifiability we have (recall that $\xi[i,j]=\xi_i\cdots\xi_j$):

\begin{lem}\cite{Mo}\label{Mosse}
There exists an integer $C$ (depending on $\sigma$) such that,  if $W\in\xi$ can be written as
$W=\xi[i-C,i+C]= \xi[j-C,j+C]$ with $i\in E_1$, then we have $j\in E_1$.
\end{lem}

We shall say that $\xi[i-C,i+C]$ and $ \xi[j-C,j+C]$ have a  
{relative common cutting position} (at the positions $i$ and $j$ respectively).
 As a consequence, if $W$ is long enough, say $|W|\ge L$  with
 \begin{equation}\label{L1}
L=\max\{2C+\max\{|A|,|B|\},|A|+|B|-1\} (>|W_0|)
\end{equation}
 and it appears at different positions in $\xi$: $W=\xi[i_1,i_2]=\xi[j_1,j_2]$, then roughly speaking, at the middle position of  $\xi[i_1,i_2]$ and $\xi[j_1,j_2]$, they have a relative common cutting position: for some integer $N \in ({|W|}/{2}-\max\{|A|,|B|\},{|W|}/{2}+\max\{|A|,|B|\})$,  $i_1+N\in E_1$ and $j_1+N \in E_1$.

\section{The Operator $T$ and Structure of  $\LS$ }

\bigskip

Define $T: \A^*\rightarrow\A^*$:  $$T(W)=W_0\sigma(W).$$

Notice  that $T$ is not a morphism on $\A^*$. It is readily checked that
$T$ is injective and
\begin{equation}\label{hhhh}
T^{n}(W)=W_0\sigma(W_0)\cdots \sigma^{n-1}(W_0)\sigma^{n}(W).
 \end{equation}

\begin{lem}\label{Tw}
   If $W\in\xi$, then $T(W)\in\xi$. Moreover, $T(W)=W_0[\sigma(W)]$.
\end{lem}
\prf Due to the primitivity of $\sigma$, the fixed sequence $\xi$ is recurrent. Thus for any $n\in\mathbb{N}$, $UW\in\xi$ for some $U\in\A^*$ with $|U|=n$. Now by the $\sigma$-invariance of $\xi$, we have that $\sigma(U)\sigma(W)\in\xi$. When the length $n$ of $U$ is large, $W_0\rt\sigma(U)$ by Lemma \ref{w0}, therefore $T(W)=W_0[\sigma(W)]\in\xi$.
\prfend

\begin{lem}\label{invert}
Let $W_1, W_2\in \A^*$. Then $T(W_1)=T(W_2)$ if and only if $W_1=W_2$; $T(W_1)\lt T(W_2)$ if and only if $W_1\lt W_2$;
$T(W_1)\rt T(W_2)$ if and only if $W_1\rt W_2$.
\end{lem}
\prf The first two easy statements hold since $\sigma$ is well marked, and the last one follows from Corollary \ref{cor21}.
\prfend

  The following lemma tells us  that if a factor $W$  appears at two positions with different natural decompositions, then,  up to a prefix $W_0'\rt W_0$, they have the same relative cutting positions.
\begin{lem}\label{2position}
Suppose that $W\in\xi$, $|W|\ge L$ with $L$ defined in (\ref{L1}), and that $W$ appears at two different positions in $\xi$, with
$W=P_1[\sigma(U_1)]Q_1$ and $W=P_2[\sigma(U_2)]Q_2$
the corresponding strict natural decompositions.
%
%
Then, denoting by $U$  the longest common suffix of $U_1$ and $U_2$ and thus writing  $U_1=U_1'U$, $U_2=U_2'U$ (where $U_1'$ or $U_2'$ is possibly empty), we have that $U$ is nonempty and
\begin{equation}\label{6600}
P_1\sigma(U_1)Q_1=W'_0[\sigma(U)]Q=P_2\sigma(U_2)Q_2,
\end{equation}
 where $Q=Q_1=Q_2$, $W_0'=P_1\sigma(U_1')=P_2\sigma(U_2')\rlt W_0$. More precisely, either $W_0'\rt W_0$, or $U_1'=U_2'=\epsilon$ and $W'_0\rt\sigma(\delta)$ for some $\delta\in\A$.

%

\end{lem}

\prf By Lemma \ref{Mosse}, the two strict natural decompositions share a
relative cutting position, and thus   all the cutting positions after this one.
This implies that $U_1$ and $U_2$ have nonempty common suffix, i.e., $U$ is
not empty. Also this implies that $Q_1=Q_2$, and consequently that
$P_1\sigma(U_1')=P_2\sigma(U_2')\rlt W_0$, where the last formula is due to Lemma \ref{w0}.
 \prfend

\begin{lem}\label{aWbW}
(1) If  $W\in\LS$ with $|W|\ge L$.
Then there exist unique $U\in \A^*, \delta\in\A$ and $Q\lt \sigma(\delta)$ with   $U\delta\in \xi$   and   $|Q|<|\sigma(\delta)|$,
 such that $$aW=aW_0[\sigma(U)]Q \quad\text{and}\quad bW=bW_0[\sigma(U)]Q.$$

  (2) If  $W\in\S$ with $|W|\ge L$.
Then there exist  $U\in \A^*$, $W_0'\in\A^*$ with either $W_0'\rt W_0$, or $W'_0\rt\sigma(\delta)$ and $|W'_0|<|\sigma(\delta)|$ for some $\delta\in\A$,
 such that $$ Wa= W_0'[\sigma(U)]a \quad\text{and}\quad  Wb= W_0'[\sigma(U)]b.$$

 (3) If $W\in \LS\cap\S$ with $|W|\ge L$. Then there exists a unique $U\in\A^*$ such that $W=T(U)$.
\end{lem}

\noindent\textbf{\emph{Remark:}}
  The word $w$ in $\LS\cap\S$ is called a bispecial word, which is developed in \cite{Ca1}, see also \cite{Ca}.

\prf (1) Consider the strict natural decompositions of $aW$ and $bW$:
$$aW=aP_a[\sigma(U_a)]Q_a \quad \text{ and}\quad bW=bP_b[\sigma(U_b)]Q_b,$$
with $U$ the longest common suffix of $U_a$ and $U_b$, $U_a=U_a'U$, $U_b=U_b'U$. Then, as in the previous proof, $U$ is nonempty, $Q_a=Q_b$,
$P_a\sigma(U_a')=P_b\sigma(U_b')$.
Moreover, putting $W_0'=P_a\sigma(U_a')$, we have that $aW_0'\rt\sigma(W_a)$ and $bW_0'\rt\sigma(W_b)$ with $W_a,W_b\in \A^*$ and the last letters of $W_a$ and $W_b$ are distinct.
Together with Lemma \ref{w0}, these facts imply that $W_0'=W_0$. %
%
%

(2) The proof for  this part is similar to the first part.

(3) This is a corollary of the first two parts.\prfend
%

\begin{lem}\label{lspecial}

(1) $W_0\in\mathcal{LS}$;

(2) Any   prefix of a left special word is left special;

(3) If $W\in\mathcal{LS}$, then $T(W)\in\mathcal{LS}$.

(4) Let $W\in\mathcal{LS}$ with $|W|\ge L$, then there exist  unique   $U\in\xi$, $\delta\in\ab$ such that
$W=W_0[\sigma(U)]Q=T(U)Q\lt T(W')$ (see Lemma \ref{aWbW}), where $W'=U\delta$. Further more, $U, W'\in\mathcal{LS}$.

\end{lem}

\prf (1) and (2) are obvious.

(3). If $aW\in\xi$, then $T(aW)\in\xi$ by Lemma \ref{Tw}.
By Lemma \ref{w0}, $\delta_aT(W)=\delta_aW_0\sigma(W)$ is a suffix of $T(aW)=W_0A\sigma(W)$, and thus $\delta_aT(W)\in\xi$.
From this, we see that $W\in\mathcal{LS}$ implies $T(W)\in\mathcal{LS}$.

(4). It follows from the  proof of the preceding lemma.\prfend

Now let $$\overline{\LS}=\bigcup\limits_{i=1}^{L}\LS_i, \quad
\overline{\LS}_n=\{W; W\lt T^n(W'), W'\in \overline{\LS}\}.$$ Remark that
$\overline{\LS}_n$ is monotone with respect to $n$.
The following theorem follows directly from the above lemma:

\bigskip
\begin{thm}\label{thm1}
$\LS=\bigcup{\overline{\LS}_n=}\lim\limits_{n\rightarrow \infty}\overline{\LS}_n$.
\end{thm}

\noindent\textbf{\emph{Remark:}} The above theorem tells us
that all left special words (which determine the complexity)
can be obtained from a finite set $\overline{\LS}$ of
left special words and by the operation $T$.

\section{Structure of $\S$ and Calculation of $\Delta^2p(n)$ }

Knowing the initial values, calculating $p(n)$ boils down into calculating $\Delta s(n+1)=\#\S_{n+1}-\#\S_{n}$. Notice that any   suffix of a special word is also special, hence if $W\in \S_{n+1}$ then $W=\delta W'$ for some $W'\in \S_{n}$ and $\delta\in \ab$. Thus the set of special words can be
 visualized as a tree showing clearly how $\S_{n+1}$ derives from $\S_{n}$
 (see the example and the figure therein in the last section).
%

As usual, for studying the special words' tree, we shall use the following notations for special words, see also \cite{CaN}:

\begin{definition}\label{02}
Let $W\in\S$.
If neither $aW$ nor $bW$ is in $\S$, we say that $W$ is a 
weak special word; If both $aW$ and $bW$ are in $\S$, we say that $W$ is a strong special word. We denote by $\S^0$ and $\S^2$ the set of weak special words and the strong weak special words respectively.
 The collection of other special words is denoted by $\S^1$.
%
%

\end{definition}

For $i\in\{0,1,2\}$,  we write $\S_n^i=\S^i\cap {\mathcal
  L}_n$. It is clear that
$$ \S_n=\S_n^0\cup\S_n^1\cup\S_n^2 {\textrm{\ and\ }}\S=\S^0\cup\S^1\cup\S^2.
$$

%
%

\begin{lem}\label{s2-s0} (1)
$\Delta s(n+1) =\  s(n+1)-s(n)=\#\S_{n}^2-\#\S_{n}^0$.

(2) $\S_n^0\cup \S_n^2 \subset \S_n\cap \LS_n.$
\end{lem}
\prf (see Theorem 4.5.4 \cite{CaN}) (1) and the fact that $\S_n^2 \subset \LS_n$ are obvious.
If a special word has only one left extension, then this left extension is also special.
\prfend

\begin{lem}\label{Twm}
  Let $c,d\in\A, W\in\xi$. If $cWd\in\xi$, then $\delta_c T(W)d\in\xi$. Conversely,
if $\delta_c T(W)d\in\xi$ and $|T(W)|\ge L$, then $cWd\in\xi$.
\end{lem}
\prf If $cWd\in\xi$, then by Lemma \ref{Tw}, $T(cWd)\in\xi$, i.e., $W_0\sigma(c)\sigma(W)\sigma(d)\in\xi$. This together with Corollary
\ref{cor21} and the fact that $\sigma$ is well marked implies that $\delta_c W_0\sigma(W)d=\delta_c T(W)d\in\xi$.

Conversely, if $\delta_c T(W)d\in\xi$ and $|T(W)|\ge L$, then by Lemma \ref{2position}, we know that $\delta_c T(W)d=\delta_c W_0[\sigma(W)]d$ is a natural decomposition. Considering the ancestor of $\delta_c T(W)d$, we know, again by Corollary \ref{cor21} and the fact that $\sigma$ is well marked, that
$cWd\in\xi$.
 \prfend

\begin{lem}\label{conserve}
If $W\in\S$, then
$T(W)\in\S$ (thus $\sigma(W)\in \S$); furthermore
$T(W)a=W_0[\sigma(W)]a$, and $T(W)b=W_0[\sigma(W)]b$.

Conversely if $W\in\S$ and $|W|\ge L$, then there exists $U\in\S$ such that $W\rt T(U)$.
\end{lem}

\prf Let $W\in\S$, then
 $Wa,Wb\in \xi$,   and by Lemma \ref{Tw},
 $$W_0[\sigma(W)]A,W_0[\sigma(W)]B\in \xi.$$
 Recalling  $A=a \circledast $ and $B=b \circledast$,
The first part of our lemma is thus proved.

The rest part is a restatement of Lemma \ref{aWbW}(2).
%
%
\prfend

%
%
%
%
%
%

We can say more on the structure of $\S^2$ and $\S^0$.

\begin{lem}\label{conserve2}
If $W\in\S^2$ then $T(W)\in\S^2$. Conversely if $W\in\S^2$ and $|W|\ge L$, then there exists a unique $U\in\S^2$ such that $W=T(U)$.
\end{lem}

\prf Let $W\in\S^2$. Then we have, by definition, that
\begin{equation}\label{123}
aWa,aWb,bWa,bWb\in\xi,
\end{equation}
 and, by Lemma \ref{Twm}, that
   $$\delta_a T(W)a, \delta_a T(W)b, \delta_b T(W)a, \delta_b T(W)b\in\xi,$$
 i.e., $T(W)\in\S^2$.  The first part of the lemma is proved.

 Now suppose $W\in\S^2$ and $|W|\ge L$. Then by Lemmas \ref{s2-s0}(2) and \ref{aWbW}(3), $W=T(U)$. By Lemma \ref{Twm}, $U\in\S^2$.
 \prfend

%


\begin{lem}\label{conserve0}
If $W\in\S^0$ and $|T(W)|\ge L$, then $T(W)\in\S^0$. Conversely if $W\in\S^0$ and $|W|\ge L$, then there exists a unique $U\in\S^0$ such that $W=T(U)$.
\end{lem}
\prf By Lemma \ref{Twm}, when $|T(W)|\ge L$ we know that $cWd\in\xi$ if and only if $\delta_c T(W)d\in\xi$. Whence $W\in\S^0$ if and only if $T(W)\in\S^0$.
The remaining proof  is almost same with the corresponding part for the preceding Lemma.
\prfend

Now denote $\overline{\S^2}=\bigcup\limits_{i=1}^{L}\S_i^2$
  the set of strong special words of length less than $L$;
  $\widetilde{\S^2}$ the set of the words $W\in \overline{\S^2} $ such that
  $|T(W)|>L$.
  The sets $\overline{\S^0}$ and $\widetilde{\S^0}$ are defined in a similar way.
  Let
\begin{equation}\label{ffff}
\widetilde{\S}=\widetilde{\S^0}\cup\widetilde{\S^2}
\end{equation}
which will be considered as ``initial special words''.

  %
%
%
%
%
\begin{lem}\label{forthm2} For any $n>L$, we have
$$\#\S_{n}^2=\sum\limits_{W\in \widetilde{\S^2}}\sum\limits_{k\ge 1}\delta(|T^k(W)|,n),
\ {\textrm{\ and\ \ \ }}
\#\S_{n}^0=\sum\limits_{W\in \widetilde{\S^0}}\sum\limits_{k\ge 1}\delta(|T^k(W)|,n),$$
where  $\delta(i,j)$ is the Kronecker symbol:
 $\delta(i,j)=1$ if $i=j$ and $=0$ otherwise.
\end{lem}
\prf Let $U\in\S_n^2$.
By Lemma \ref{conserve2},   there exist $k\ge 1$ and $W\in \widetilde{\S^2}$, which are unique, such that
$U=T^{k}(W)$. Conversely if $|T^{k}(W)|=n$ for some $k\ge 1, W\in \widetilde{\S^2}$, then
$T^{k}(W)\in\S_n^2$. Thus we have
$$
\S_n^2=\{U; U=T^{k}(W),|T^{k}(W)|=n, k\ge 1, W\in \widetilde{\S^2} \}
$$
where $k$ and $W$ in the representation  $U=T^{k}(W)$ are uniquely determined by $U$. The first equality is thus proved. The second is proved similarly. \prfend

The following formula then follows from the above lemma and Lemma \ref{s2-s0}:

\begin{lem}\label{thm2} For any $n>L$, we have
\begin{equation*}
  \begin{split}
      \Delta s(n+1) =&~s(n+1)-s(n)\\
 =& \sum\limits_{W\in\ \widetilde{\S^2}}\sum\limits_{k\ge 1}\delta(|T^k(W)|,n)
-\sum\limits_{W\in\ \widetilde{\S^0}}\sum\limits_{k\ge 1}\delta(|T^k(W)|,n).\\
  \end{split}
\end{equation*}
 It can be written as
$$\Delta s(n+1)=\sum\limits_{W\in\  \overline{\S}}
\sum\limits_{k\ge 1}\sgn(W)\delta(|T^k(W)|,n),$$
where $\overline{\S}=\bigcup\limits_{i=1}^{L}\S_i$ (the special words of length less than $L$), and
\begin{equation}\label{signw}
\sgn(W)=\left\{
\begin{array} {rl}
   -1 &  \quad  \text{ if ~ } W \in\widetilde{\S^0}   \\
    1 &  \quad  \text{ if ~ } W \in\widetilde{\S^2}  \\
    0 &  \quad  \text{ otherwise. }
\end{array}
\right.
\end{equation}

\end{lem}

\bigskip
\noindent\textbf{\emph{Remark:}}
1. The function $\sgn(\cdot)$ is equal to the bilateral multiplicity of a factor (\cite{CaN}). See Theorem 4.5.4 \cite{CaN} for more general cases.

2.
The above lemma tells us that the complexity $p(n)$ can be computed knowing a finite set $\overline{\S}$ of special words. In the next section, we will find out a (non-linear) recurrence formula for the computation.

\section{Recurrence Formula for the Complexity}

\bigskip

Recall that $M$ denotes the incidence matrix of $\sigma$. Then $M^2$ is the incidence matrix of $\sigma^2$ which possess non-negative eigenvalues.  Since $\sigma$ and $\sigma^2$ share the   fixed sequence $\xi$,  we may suppose  without  loss of generality that
$\textrm{the eigenvalues of}\  M\ \textrm{is non-negative.}$

Let $\lambda_1\ge\lambda_2\ge0$ be the two eigenvalues, $V_1$, $V_2$ be the corresponding eigenvectors. Since $M$ is primitive, $\lambda_1>\lambda_2$ and  $V_1$ is positive.

Recall that: for $W\in\ab^*$, $L(W)=(|W|_a,|W|_b)^t$,
\begin{equation}\label{gggg}
|\sigma^n(W)|=(1,1)M^nL(W).
\end{equation}

\begin{lem}\label{sync}
Let $X,Y\in\R^2$. Then there exists $N=N(X,Y)\ge1$ such that $(1,1)M^{N+n}(X-Y)$ ($n\in\mathbb{N}$) is of constant sign. That is,
$$(1,1)M^{N+n}X> (\textrm{resp.\ } =,<) \ (1,1)M^{N+n}Y  \textrm{~for all\ } \ n\in\mathbb{N} .$$

\end{lem}

\prf Let $X-Y=\mu_1V_1+\mu_2V_2$ where $\mu_1, \mu_2\in \R$, then for $k\ge 1$,
$$(1,1)M^k(X-Y)=\lambda_1^k\mu_1(1,1)V_1+\lambda_2^k\mu_2(1,1)V_2.
$$

\smallskip

\noindent Case 1. $\mu_1=0$. Then $(1,1)M^k(X-Y)=\lambda_2^k\mu_2(1,1)V_2$, which is obviously of the sign of $\lambda_2\mu_2(1,1)V_2$ independent of $k\ge1$.

\smallskip

\noindent Case 2. $\mu_1>0$. Since $\lambda_1>0$, $(1,1)V_1>0$ and $\lambda_1>\lambda_2\ge0$, there exists $N\ge 1$ such that for $k\ge N$ we have
$\lambda_1^k\mu_1(1,1)V_1+\lambda_2^k\mu_2(1,1)V_2>0
$.

\smallskip

\noindent Case 3.   $\mu_1<0$. The similar proof as Case 2.
\prfend

\begin{cor}\label{cor1}  Let $W_1,W_2\in\A^*$. There exists $N=N(W_1,W_2)$ such that
$|T^{N+n}(W_1)|-|T^{N+n}(W_2)|\  (n\in\mathbb{N})$ is of constant sign. This sign (called the final sign) will be denoted by $\fsgn\{W_1,W_2\}$.

\end{cor}
\prf   The lemma follows directly from the above lemma and (\ref{gggg}).
 \prfend


In fact, we can say more:

\begin{cor}\label{T1} Let $W_1,W_2\in\A^*$. Then there exist $m_1=m_1(W_1,W_2), m_2=m_2(W_1,W_2)\in \mathbb{N}$ such that one of the following alternatives holds:

 (1).  $|T^{m_1}(W_1)|=|T^{m_2}(W_2)|<|T^{m_1+1}(W_1)|=|T^{m_2+1}(W_2)|<|T^{m_1+2}(W_1)|=|T^{m_2+2}(W_2)|<\ldots$

(2).  $|T^{m_1}(W_1)|<|T^{m_2}(W_2)|<|T^{m_1+1}(W_1)|<|T^{m_2+1}(W_2)|<|T^{m_1+2}(W_1)|<|T^{m_2+2}(W_2)|<\ldots$.

 \end{cor}

\prf   If $\fsgn(T^{m}(W_1), T^{n}(W_2))=0$ for some $m,n \in \mathbb{N}$, the alternative (1) holds.

     Otherwise,  $\fsgn(T^{m}(W_1), T^{n}(W_2))\not=0$ for any $m, n\in \mathbb{N}$.
      We assume, without loss of generality, that $ \fsgn(W_1, W_2)=-1$.
      Due to the primitivity, $W_2$ is a factor of $T^l(W_1)$ for $l$
      large enough, and it turns out that  $\fsgn(T^{l}(W_1), W_2)=1$.
       Now clearly $m\mapsto \fsgn(T^m(W_1), W_2)$ is an increasing mapping
       from $\mathbb{N}$ onto $\{-1,1\}$,
     therefore there exists $m\in \mathbb{N}$ 
      such that $\fsgn(T^m(W_1), W_2)=-1$, while $\fsgn(T^{m+1}(W_1), W_2)=1$. Whence the alternative (2) holds for $$m_2=\max\{N(T^m(W_1), W_2), N(T^{m+1}(W_1), W_2)\}, {\rm  and ~ } m_1=m+m_2.$$
\vskip -0.6cm
\prfend

Now we can deduce from the above lemma the recurrence properties of the complexity. First let $\widetilde{\S}=\{S_1,S_2,\cdots,S_K\}$    and denote
$$n_1-1=\max\big\{\max\{m_1(W_1,W_2), m_2(W_1,W_2)\};\ \  W_1,W_2\in\widetilde{\S}\big\},$$
where $m_1(W_1,W_2), m_2(W_1,W_2)$ are defined in Lemma \ref{T1}.

 We start from $T^{n_1}(S_1)$. By Lemma \ref{T1}, for each $j=2,3,\cdots,K$, there exists unique $n_j\in \mathbb{N}$ such that $|T^{n_1}(S_1)|\le |T^{n_j}(S_j)|<|T^{n_1+1}(S_1)|$. Without loss of generality we may suppose that
$$|T^{n_1}(S_1)|\le |T^{n_2}(S_2)|\le \cdots \le |T^{n_K}(S_K)\le |T^{n_1+1}(S_1)|. $$

Then  for simplifying the notations let  $N_k^j=|T^j(T^{n_k}(S_k))| \ \ (1\le k\le K,\ j\in \mathbb{N})$.
We have  by  Lemma \ref{T1} the following unison property for the ``jumps of $|T^i(W_k)|$":
\begin{eqnarray}\label{last}
      & N_1^0&\le N_2^0\cdots \le N_K^0  \nonumber \\
  \le & N_1^1&\le N_2^1\cdots \le N_K^1   \nonumber \\
      &\cdots  &   \cdots                      \nonumber      \\
 \le  & N_1^j&\le N_2^j\cdots \le N_K^j   \nonumber  \\
 \le  & N_1^{j+1}&\le\cdots
\end{eqnarray}

Now we can formulate the recurrence formula of the complexity.
Let $\chi_{[m,n)}$ denote the indicator function of the integers' interval $[m,n)$.
Let $I^j=[N_1^j, N_1^{j+1}), j\in \mathbb{N}$. We see that $I^j $ is the disjoint union of the  subintervals     $I_k^j=[ N_k^j, N_{k+1}^{j})\  (j\in \mathbb{N}, k\in\{1,2,\cdots,K\})$, where  $N_{K+1}^{j}=N_{1}(j+1)$. That is
$$[N_1^0,\infty)=\bigcup\limits_{j=0}^{\infty} I^j,\ \ I^j=\bigcup\limits_{k=1}^K I_k^j\ .$$

 \medskip
\subsection{Initial values of the complexity}

\mbox{}

\smallskip

 Finally let $c_k=\sum\limits_{i=1}^{K} \sgn(S_i)\delta(|T^{n_i}(S_i)|,|T^{n_k}(S_k)|)$ ($k=1,\cdots,K$), where $\sgn(\cdot)$ is defined in (\ref{signw}).
Then by Lemma \ref{thm2}, we have,
$\Delta s(n+1)=c_k$   if $n=|T^{n_k}(S_k)| (k=1,\cdots,K)$ and $=0$ otherwise.  In other words, $n\mapsto s(n+1)\ (n\in I^0)$   is a step function with jumps $c_k$ at $n=N_k(0)$\
 $(k\in\{1,2,\cdots,K\})$:
\begin{equation}\label{unison2}
s(n+1)=s(N_1(0))+\sum\limits_{k=1}^K(c_1+\cdots+c_k)\chi_{I_k^0}(n) \ (n\in I^0),
\end{equation}

\medskip
\subsection{Recurrence formula of $s(\cdot +1)$ on $I^j$} \mbox{}

\smallskip

 Notice that $I^j=\bigcup\limits_{k=1}^K I_k^j$ ($j\in\mathbb{N}$) can be  calculated  directly or by some easy  recurrence formula as described in the following:

\begin{prp}\label{intv} We have for any $W\in\A^*$, $n\in\mathbb{N}$,

1. $|\sigma^{n+2}(W)|=\tr(M)\ |\sigma^{n+1}(W)|-\det(M)\ |\sigma^{n}(W)|$;

\quad  \! $|T^{n+2}(W)|=\tr(M)\ |T^{n+1}(W)|-\det(M)\ |T^{n}(W)|+a,$

where $a=|\sigma(W_0)|-(\tr(M)-1)|W_0|$.

2.
$|\sigma^{n}(W)|=\lambda_1^n\mu_1(1,1)V_1+\lambda_2^n\mu_2(1,1)V_2$ if $L(W)=\mu_1V_1+\mu_2V_2$
;

\quad  \! $|T^{n}(W)|=\lambda_1^n\mu_1(1,1)V_1+\lambda_2^n\mu_2(1,1)V_2+b_n$,

 where $b_n (n\in\mathbb{N})$ is a fixed sequence  given explicitly by $L(W_0)$ and $M$.
\end{prp}

 \prf All the results can be deduced easily from (\ref{hhhh}), (\ref{gggg}) and  Cayley-Hamilton formula (with $I$ denotes the identity matrix): $M^2=\tr(M)\ M-\det(M)\ I$.
 \prfend

 \bigskip
We have just seen the recurrence properties of the intervals $I_j$ ($j\in\mathbb{N}$). Still using  Lemma \ref{thm2} and the formula (\ref{last}) and we see that what happens for $s(n+1)$ ($n\in I^j$, $j\in\mathbb{N}$) is recurrently the same as  $s(n+1)$ ($n\in I^0$), i.e., similar to (\ref{unison2}) we have proved the following 

\bigskip
\begin{thm}\label{thm3} Let $\sigma$ be a well marked, primitive, non-periodic substitution having non-negative eigenvalues. Then for $n\in [N_1^0,\infty)=\bigcup\limits_{j=0}^{\infty} I^j$,  the following recurrence formula holds:
$$s(n+1)=s(N_1^j)+\sum\limits_{k=1}^K(c_1+\cdots+c_k)\chi_{I_k^j}(n) \  (n\in I^j, j\ge0).$$
\end{thm}

\bigskip
\noindent\textbf{\emph{Remark:}}
1. The conditions ``primitive, well marked, non-periodic, having non-negative eigenvalues'' are non-essential as have already mentioned.


2. $s(N_1^{j+1})-s(N_1^j)\equiv c_1+\cdots+c_K $ ($j\in\mathbb{N}$), which implies roughly $s(\lambda_1^n )\approx n(c_1+\cdots+c_K) $ for large $n$.

3. Although the above mentioned $N_1^0$ can be more or less controlled in the proof of the theorem, but how to give efficiently this big integer $N$ remains as an open problem.

\bigskip
Fially let us give briefly an example:
consider the substitution $\sigma=(aab,ba)$ i.e., $a\mapsto aab, b\mapsto ba$.

\begin{figure}
\begin{tikzpicture}[scale=0.95]
\tikzset{grow=down,level distance=26pt}
\Tree[.$\varepsilon$
 [.$a$ [.$aa$ [.$baa$ $abaa$  ]
       ]
       [.$ba$ [.$bba$ [.$abba$ [.$aabba$ [.$baabba$ [.$abaabba$ [.$aabaabba$  [.$aaabaabba$ \mbox{}  ] [.$baabaabba$ \mbox{} ]   ] ] ] ]
                      ]
              ]
       ]
 ]
 [.$b$ [.$ab$ [.$aab$ [ .$aaab$ [.$baaab$ [.$abaaab$ [.$babaaab$ [.$bbabaaab$ [.$abbabaaab$ \mbox{} ]] ] ]]]
                      [ .$baab$ [.$abaab$ [.$aabaab$ [.$aaabaab$ [.$baaabaab$ [ .$bbaaabaab$ \mbox{} ] ] ] ]]]
       ] ]
 ]
]
\end{tikzpicture}
  \caption{Tree of Special Words}\label{sp}
\end{figure}
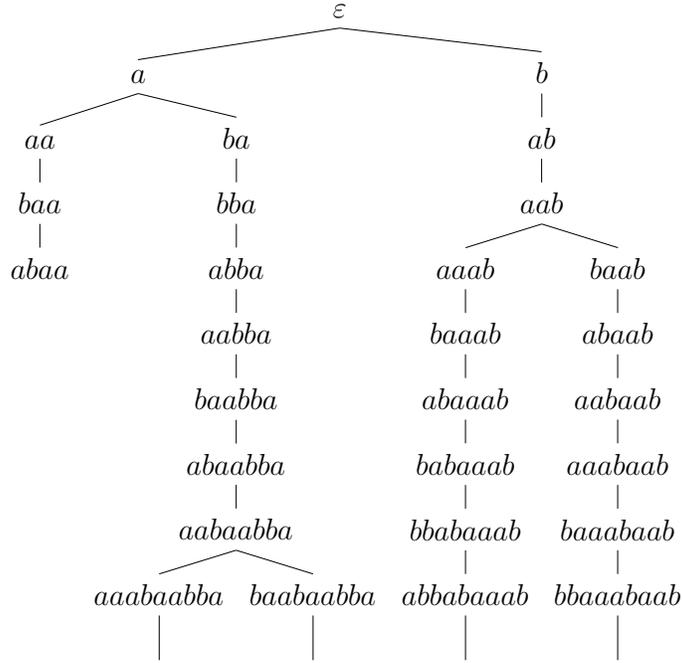

For this substitution, we have $W_0=\varepsilon$ and thus $T=\sigma$. The incidence matrix $M=\left(
                         \begin{array}{cc}
                           2 & 1 \\
                           1 & 1 \\
                         \end{array}
                       \right)
$ and the characteristic polynomial is $\lambda^2-3\lambda+1$.
The fixed point reads $$\xi=aabaabbaaabaabbabaaabaabaabbaaabaabbabaaabbaaabaab\cdots$$
The tree of the special words is depicted in Figure \ref{sp}.

The weak and strong special words (here $\sigma^0$ is the identity map):

\quad $\S^0=\{abaa,aabbaaabaab,\cdots\}=\{\sigma^n(abaa); n=0,1,2,\cdots\}$,

\quad $\S^2=\{\varepsilon, a,aab,aabaabba,\cdots\}=\{\varepsilon\}\cup\{\sigma^n(a); n=0,1,2,\cdots\}$.

From the structure of special words, the numbers of special words $s(n)$ and the complexity $p(n)$ read

\begin{center}
\begin{tabular}{c|c|c|c|c|c|c|c|c|c|c|c|c|c|c|c|c|c}
\hline
  $n$ &0 &1 & 2 & 3 & 4 & 5 &6 &7&8&9&10&11&12&13&14&15&$\cdots$   \\
  \hline
  $s(n)$&1 & 2 & 3 & 3 & 4 & 3 &3 &3&3&4&4&4&3&3&3&3&$\cdots$   \\
\hline
$p(n)$&1 &2 & 4 &7 & 10 & 14 &17 &20&23&26&30&34&38&41&44&47&$\cdots$   \\
\hline
\end{tabular}
\end{center}

We can formulate $s(n)$ as
$$s(n)=\begin{cases}
   1 & \text{if } n=0,\\
   2 & \text{if } n=1,\\
   3 & \text{if } n\in\{2,3 \}\cup \bigcup_{k\ge0} [d(k)+1,g(k+1)],\\
   4 & \text{if } n\in\bigcup_{k\ge0} [g(k)+1,d(k)],
   \end{cases}
$$
where the number sequences $g(k)$ and $d(k)$ are defined as
$$g(k)=(1,1)M^k(2,1)^t,\quad d(k)=(1,1)M^k(3,1)^t,$$
satisfying both the same recurrence:
$$\begin{cases}
 ~   g(k+2)=3g(k+1)-g(k), \\ ~d(k+2)=3d(k+1)-d(k),\\
\end{cases}$$
with $g(0)=3,g(1)=8$ and $d(0)=4,d(1)=11.$

\bigskip

\noindent \textsc{Acknowlegement}\quad
The authors would like to thank Prof. Z.Y. Wen (Tsinghua),\ J.P. Allouche (Jussieu) and others for helpful discussions,references and corrections.

\end{document}